\documentclass[5p,fleqn]{elsarticle}
\usepackage{graphicx}
\usepackage{mathtools}
\usepackage{amssymb}

\usepackage{listings}
\usepackage{xcolor}
\usepackage{widetext}
\usepackage{cuted}

\newcommand*\diff{\mathop{}\!\mathrm{d}}
\newcommand{\CC}{C\nolinebreak\hspace{-.05em}\raisebox{.4ex}{\tiny\bf +}\nolinebreak\hspace{-.10em}\raisebox{.4ex}{\tiny\bf +} }

\definecolor{commentsColor}{rgb}{0.497495, 0.497587, 0.497464}
\definecolor{keywordsColor}{rgb}{0.000000, 0.000000, 0.635294}
\definecolor{stringColor}{rgb}{0.558215, 0.000000, 0.135316}

\allowdisplaybreaks

\begin{document}
\title{A simple algorithm for uniform sampling on the surface of a hypersphere}

\author{Stefan Schnabel}
\ead{stefan.schnabel@itp.uni-leipzig.de}

\author{Wolfhard Janke}
\ead{wolfhard.janke@itp.uni-leipzig.de}

\address{Institut f\"ur Theoretische Physik, Universit\"at Leipzig, IPF 231101, 04081 Leipzig, Germany}
\date{\today}





\begin{abstract}

We propose a simple method for uniform sampling of points on the surface of a hypersphere in arbitrarily many dimensions. By avoiding the evaluation of computationally expensive functions like logarithms, sines, cosines, or higher-order roots the new method is faster than alternative techniques.

\end{abstract}

\maketitle

\section{Introduction}

The (pseudo) random sampling of points uniformly distributed on the surface of an $n$-dimensional sphere is a problem that is regularly encountered in the context of probalistic algorithms. These can, for instance, be Monte Carlo simulations in the context of statistical physics, evolutionary algorithms in high-dimensional parameter spaces, or methods to generate random rotation matrices. The easiest and most efficient way to solve this problem that is used nowadays was suggested by Muller \cite{Muller} and prescribes to draw all coordinates independently from a normal distribution exploiting that the multivariate Gaussian distribution is spherically symmetric. Then, the resulting vector only has to be rescaled such that its length equals the radius of the sphere.

Although other techniques, e.g., \cite{Hicks,Sibuya,Guralink,Tashiro,Yang}, have been introduced, they are more complicated or do not appear to improve the performance such that they do not enjoy much popularity. The case is altered for low number of dimensions, where on the one hand the so-called rejection method becomes viable. It generates vectors to random points in the hypercube until a point is obtained that happens to lie within the hypersphere thus creating a spherically symmetric distribution. The respective vector is then rescaled to the desired length. On the other hand for $n=3$ the polar method is widely used and Marsaglia has introduced very effective algorithms for $n=3$ and $n=4$ \cite{Marsaglia} all of which we will discuss below. Since so far, no generalizations of Marsaglia's methods are publicly known \cite{Wolfram,NumRec} we will show here how a very similar method can be developed for an arbitrary number of dimensions.

\subsection{Polar method}
In three dimensions it is a well-known fact that any single cartesian coordinate of points on the sphere is uniformly distributed. This is used for what can be called the polar method:
\begin{itemize}
\item{Draw $X_3\in[-1,1)$ and $\phi\in[0,2\pi)$.}
\item{Assign remaining coordinates:
    \begin{eqnarray}
        X_1&\coloneqq&\cos(2\pi\phi)\sqrt{1-X_3^2},\\
        X_2&\coloneqq&\sin(2\pi\phi)\sqrt{1-X_3^2}.
    \end{eqnarray}
}
\end{itemize}

\subsection{Marsaglia $n=3$}
However, the use of the comparatively slow trigonometric functions can be avoided \cite{Marsaglia}. Using the rejection method to generate a point within a (two dimensional) circle one obtains a a random direction while the distribution of the distance $r$ from the origin is proportional to $r$ itself. This means that $r^2$ is uniformly distributed in $[0,1)$ and can, therefore, be used to define the third coordinate.
\begin{itemize}
\item{Draw $(a,b)\in[-1,1)^2$ until $1>S\coloneqq r^2= a^2+b^2$.}
\item{Assign coordinates:
    \begin{eqnarray}
        X_3&\coloneqq& 1-2S,\\
        X_1&\coloneqq& a\sqrt{\frac{1-X_3^2}{S}},\\
        X_2&\coloneqq& b\sqrt{\frac{1-X_3^2}{S}}.
    \end{eqnarray}
}
\end{itemize}

\subsection{Marsaglia $n=4$}
Using random points within a circle is even more natural in four dimensions, since there, as can easily be shown, any two cartesian coordinates follow this distribution anyway.
\begin{itemize}
\item{Draw $(X_1,X_2)\in[-1,1)^2$ directly until $1>S_1\coloneqq X_1^2+X_2^2$.}
\item{Draw $(a,b)\in[-1,1)^2$ until $1>S_2\coloneqq a^2+b^2$.}
\item{Assign coordinates $X_3$ and $X_4$:
    \begin{eqnarray}
        X_3&\coloneqq& a\sqrt{\frac{1-S_1}{S_2}},\\
        X_4&\coloneqq& b\sqrt{\frac{1-S_1}{S_2}}.
    \end{eqnarray}
}
\end{itemize}

\section{Our method}

Our algorithm to create a random point $\mathbf{X}$ on an $n$-sphere with even $n$ consists of the following three simple steps:
\begin{itemize}
\item{Draw $n/2$ pairs of coordinates $(a,b)\in[-1,1)^2$ until $1>a^2+b^2$.}
\item{Sort the pairs such that $S_{i-1}\le S_{i}$ with $S_i=a_i^2+b_i^2$ and set $S_0\coloneqq0$.}
\item{Assign coordinates $X$:
    \begin{eqnarray}
        X_{2i-1}&\coloneqq& r_ia_i, \label{eqn:oM_Xi_1}  \\
        X_{2i}  &\coloneqq& r_ib_i,
    \end{eqnarray}
    for $i=1,\dots,n/2$ where
    \begin{equation}
    	r_i=\sqrt{\frac{1-S_{i-1}/S_{i}}{S_{n/2}}}.
    	\label{eqn:oM_ri_S} 
    \end{equation}
}
\end{itemize}
The method is computationally very efficient for small or moderately large $n$ since no trigonometric or logarithmic functions and no roots of higher-order are necessary.
For the case of odd $n$ we follow the suggestion in \cite{Sibuya} of creating a random vector in $n+1$ dimensions and normalizing its first $n$ components.

A useful variant of the method is obtained if (\ref{eqn:oM_ri_S}) is replaced by
\begin{equation}
    	r_i=\sqrt{1-\frac{S_{i-1}}{S_{i}}}.
    	\label{eqn:ri_B} 
\end{equation}
Then $\mathbf{X}$ becomes a point that is randomly distributed in the interior of the $n$-sphere, i.e., in the $n$-ball. Unfortunately, we do not see how a random point in the $n$-ball with odd $n$ can be obtained in a similar way. A well-known alternative that works for general $n$ is to create a random point on the $(n+2)$-sphere and use the first $n$ coordinates of that point.

Although we developed the method from first principles, in its final form it appears to be closely related to Sibuya's method \cite{Sibuya} and it is, therefore, fitting that we first show how our method can be derived from the latter.

\subsection{Sibuya's method (for even $n$)}
\begin{itemize}
	\item{Obtain $n/2$ points $(x_i,y_i)$ randomly distributed \emph{on} a unit-circle, either by drawing an angle or using the rejection technique.}
	\item{Draw $n/2-1$ random variables $U\in[0,1)$ and sort them:
	    \begin{equation}
	    	0=U_0\le U_1\le \dots \le U_{n/2-1}\le U_{n/2}=1.
	    \end{equation}
	}
	\item{Assign coordinates $X$:
	    \begin{equation}
	        (X_{2i-1},X_{2i}) \coloneqq  (x_i,y_i)\sqrt{U_i-U_{i-1}}, \label{eq:sibuya_coord_1}\\
	    \end{equation}
	    for $i=1,\dots,n/2$.
	}
\end{itemize}
\vspace{.0cm}

\subsection{Modifying Sibuya's method}
We first change the method slightly, by using $n/2$ random variables $U$ such that $U_{n/2}=\max(\{U\})<1$. Consequently, by rescaling the differences $U_i-U_{i-1}$ the factor in (\ref{eq:sibuya_coord_1}) becomes $\sqrt{(U_i-U_{i-1})/U_{n/2}}$.

We exploit the important yet simple observation mentioned above: If the rejection method is used, i.e., if a random point in a square $(a,b)\in[-1,1)^2$ is repeatedly drawn until $a^2+b^2<1$ and then normalized $(x,y)=(a,b)/\sqrt{a^2+b^2}$, then the squared radius 
$a^2+b^2$ is uniformly distributed in $[0,1)$ and independent of $x$ and $y$. Therefore, if we use the rejection method to determine the points $(x_i,y_i)=(a_i/\sqrt{a_i^2+b_i^2},b_i/\sqrt{a_i^2+b_i^2})$ we can also use the radii to obtain a sequence of random numbers $V_i=a_i^2+b_i^2$ which is sorted to create the sequence $U_i$.

Equation (\ref{eq:sibuya_coord_1}) can now be written as  
    \begin{eqnarray}
        \!\!(X_{2i-1},X_{2i})&\coloneqq& (x_i,y_i)\sqrt{\frac{U_i-U_{i-1}}{U_{n/2}}},\nonumber\\
          &=&\frac{(a_i,b_i)}{\sqrt{V_{i}}} \sqrt{U_{i}} \sqrt{\frac{1-U_{i-1}/U_{i}}{U_{n/2}}}.
    \end{eqnarray}
Since the numbers $\{U_i\}$ and the coordinates $\{(x_i,y_i)\}$ are independent, we can arbitrarily permutate the latter and in particular arrange them in such a way that each pair of coordinates is again matched with its original radius. Then $\sqrt{U_{i}}$ cancels with $\sqrt{V_{i}}$ and we obtain our method (\ref{eqn:oM_Xi_1})-(\ref{eqn:oM_ri_S}) as prescribed. We thus improve on Sibuya's original method by avoiding the independent drawing of the random numbers $U$ and by only requiring half as many square roots.

\subsection{Our method from scratch}
It may be instructive to see how we originally derived the method from basic principles. In the following $\mathcal{S}_{n}$ denotes an $n$-sphere while $\mathcal{B}_{n}$ stands for its interior, i.e., the unit ball.

Consider the probability density of points in $\mathcal{S}_{n}$:
\begin{equation}
P_n(X_1,\dots,X_n)=\frac1{A_n}\delta(X_1^2+\dots+X_n^2-1),
\end{equation}
where $A_n$ is a normalizing factor. Integrating over two coordinates,
\begin{widetext}
    \begin{eqnarray}
        \int\limits_{-\infty}^{\infty}\int\limits_{-\infty}^{\infty} P_n(\dots)\diff X_{n-1}\diff X_{n} &=&\frac1{A_n}\int\limits_{-\infty}^{\infty}\int\limits_{-\infty}^{\infty}\delta(X_1^2+\dots+X_n^2-1)\diff X_{n-1}\diff X_{n}, \nonumber\\
        &=&\frac1{A_n}\int\limits_{0}^{\infty}\int\limits_{0}^{2\pi}\delta(X_1^2+\dots+X_{n-2}^2+r^2-1)r\diff \phi\diff r, \nonumber\\
        &=&\frac{\pi}{A_n}\int\limits_{0}^{\infty}\delta(X_1^2+\dots+X_{n-2}^2+\rho-1)\diff \rho, \nonumber\\
        &=& 
            \begin{cases}
                \frac{\pi}{A_n} & \text{if} \quad X_1^2+\dots+X_{n-2}^2<1 \\
                0                   & \text{if} \quad X_1^2+\dots+X_{n-2}^2>1
            \end{cases},
        \label{eq:unif_dis}
    \end{eqnarray}
\end{widetext}
demonstrates the known fact that if points uniformly distributed in $\mathcal{S}_{n}$ are projected onto a $(n-2)$-plane their images are uniformly distributed in $\mathcal{B}_{n-2}$. It is also clear due to symmetry that for the subset of points in $\mathcal{S}_{n}$ with coordinates $X_1,\dots,X_{n-2}$ set to some fixed values the remaining two coordinates are uniformly distributed on a ring in the $X_{n-1},X_n$-plane of radius $(1-X_1^2+\dots+X_{n-2}^2)^{1/2}$ and can, therefore, easily be sampled with the correct distribution. If we assume that a method of generating a random point in $\mathcal{S}_{n-2}$ is available, a random point in $\mathcal{B}_{n-2}$ is obtained by rescaling to a random radius $r_{n-2}\in[0,1)$ with a probability density $P_{r_{n-2}}(x)\propto x^{n-3}$. Hence, if $(Y_1,\dots,Y_{n-2})$ is a random point in $\mathcal{S}_{n-2}$ and $(x,y)$ a random point in $\mathcal{S}_{2}(1)$, i.e., on the unit circle, then $(r_{n-2}Y_1,\dots,r_{n-2}Y_{n-2},x(1-r_{n-2}^2)^{1/2},y(1-r_{n-2}^2)^{1/2})$ is a random point in $\mathcal{S}_{n}$.

For even $n$ we thus obtain by recursion:
\begin{eqnarray}
(X_1,X_2) & = & (x_1,y_1) \cdot r_2 \cdot r_4 \cdot r_6 \dots r_{n-2}, \nonumber\\
(X_3,X_4) & = & (x_2,y_2) \sqrt{1-r_2^2} \cdot r_4 \cdot r_6 \dots r_{n-2}, \nonumber\\
(X_5,X_6) & = & (x_3,y_3) \sqrt{1-r_4^2} \cdot r_6 \dots r_{n-2}, \nonumber\\
    & \vdots & \\
(X_{n-1},X_n) & = & (x_{n/2},y_{n/2}) \sqrt{1-r_{n-2}^2}, \nonumber
\end{eqnarray}
This leaves us with the simple task of generating points $(x_i,y_i)\in\mathcal{S}_{2}$ and the less trivial question of how to obtain the radii $r_{2i}$. Since their probability densities follow simple power laws, they could be obtained from uniformly distributed random variables $\xi_i\in[0,1)$ through $r_{2i}=\xi_i^{1/2i}$, yet this is not suitable due to the high computational effort required to calculate higher-order roots. Another way to obtain random variables with a power-law distribution is to consider the maximum of multiple independent random variables.

If $m$ random variables $\xi_i\in[0,1)$ are independent and uniformly distributed: $P_\xi(x)= \text{const}$ then the probability of their maximum $P_{\max(\xi_1,\dots,\xi_m)}(x)\propto x^{m-1}$. Or more general, if $m$ random variables $\xi_i\in[0,1)$ are independent and distributed according to a power law, $P_\xi(x)\propto x^\alpha$, then the probability of their maximum follows $P_{\max(\xi_1,\dots,\xi_m)}(x)\propto x^{(\alpha+1)m-1}$.

It follows that if $m$ independent random variables $\xi\in[0,1)$ are sorted, such that $\xi_1\le\xi_2\le\dots\le\xi_m$, then the random variables $\eta_i=\xi_i/\xi_{i+1}\in[0,1]$ are independent and are distributed according to $P_{\eta_i}(x)\propto x^{i-1}$. Furthermore, if $m$ independent random variables $\xi\in[0,1)$ which are distributed according to a power law, $P_\xi(x)\propto x^\alpha$, are sorted, then the random variables $\eta_i=\xi_i/\xi_{i+1}\in[0,1]$ are independent and are distributed according to $P_{\eta_i}(x)\propto x^{(\alpha+1)i-1}$.

In our case, this implies that we require $n/2$ random variables $\xi\in[0,1)$ with  $P_\xi(x)\propto x$, i.e., $\alpha=1$. These are to be sorted such that $\xi_1\le\xi_2\le\dots\le\xi_{n/2}$ and the radii are obtained through $r_{2i}=\xi_i/\xi_{i+1}$. For the resulting vector we obtain
    \begin{eqnarray}
    (X_1,X_2) & = & (x_1,y_1) \cdot \frac{\xi_1}{\xi_2} \cdot \frac{\xi_2}{\xi_3} \cdot \frac{\xi_3}{\xi_4} \dots \frac{\xi_{n/2-1}}{\xi_{n/2}},\nonumber\\
    (X_3,X_4) & = & (x_2,y_2) \sqrt{1-\left(\frac{\xi_1}{\xi_2}\right)^2} \cdot \frac{\xi_2}{\xi_3} \cdot \frac{\xi_3}{\xi_4} \dots \frac{\xi_{n/2-1}}{\xi_{n/2}},\nonumber\\
    (X_5,X_6) & = & (x_3,y_3) \sqrt{1-\left(\frac{\xi_2}{\xi_3}\right)^2} \cdot \frac{\xi_3}{\xi_4} \dots \frac{\xi_{n/2-1}}{\xi_{n/2}},\nonumber\\
        & \vdots & \\
    (X_{n-1},X_n) & = & (x_{n/2},y_{n/2}) \sqrt{1-\left(\frac{\xi_{n/2-1}}{\xi_{n/2}}\right)^2},\nonumber
    \end{eqnarray}
which simplifies to
    \begin{eqnarray}
    (X_1,X_2) & = & (x_1,y_1) \cdot \frac{\xi_1}{\xi_{n/2}},\nonumber\\
    (X_3,X_4) & = & (x_2,y_2) \sqrt{1-\left(\frac{\xi_1}{\xi_2}\right)^2} \cdot \frac{\xi_2}{\xi_{n/2}},\nonumber\\
    (X_5,X_6) & = & (x_3,y_3) \sqrt{1-\left(\frac{\xi_2}{\xi_3}\right)^2} \cdot \frac{\xi_3}{\xi_{n/2}},\nonumber\\
        & \vdots & \\
    (X_{n-1},X_n) & = & (x_{n/2},y_{n/2}) \sqrt{1-\left(\frac{\xi_{n/2-1}}{\xi_{n/2}}\right)^2} \cdot \frac{\xi_{n/2}}{\xi_{n/2}},\nonumber
    \end{eqnarray}
By defining $\xi_0\coloneqq0$ we can use the general expression
\begin{equation}
\begin{array}{lcl} 
(X_{2i-1},X_{2i}) & = & (x_i,y_i) \sqrt{1-\left(\frac{\xi_{i-1}}{\xi_i}\right)^2} \cdot \frac{\xi_i}{\xi_{n/2}},
\end{array}
\label{eqn:deriv_scr_a}
\end{equation}
for $i=1,\dots,n/2$.

If we use the rejection method to determine the points $(x_i,y_i)\in\mathcal{S}_{2}$, i.e., if we repeatedly draw $(a_i,b_i)\in[-1,1)^2$ until $a_i^2+b_i^2<1$ and assign $x_i\coloneqq a_i/\sqrt{a_i^2+b_i^2}, y_i\coloneqq b_i/\sqrt{a_i^2+b_i^2}$ then, as previously discussed, the radii $R_i=\sqrt{a_i^2+b_i^2}$ are independent of $(x_i,y_i)$ and, since they have linearly growing probability distribution, $P_R(x)\propto x$, are a suitable choice for $\{\xi_i\}$ once they are sorted.

Finally, we note that the order in which the points $(x_i,y_i)$ appear in (\ref{eqn:deriv_scr_a}) is irrelevant and can be changed at will. In particular, we can permutate the points such that each one is again matched up with its original radius. Hence, if $(\hat{x}_i,\hat{y}_i) \coloneqq (\hat{a}_i,\hat{b}_i)/\sqrt{\hat{a}_i^2+\hat{b}_i^2}$ are the sorted sequences $\hat{a}_1^2+\hat{b}_1^2<\dots<\hat{a}_{n/2}^2+\hat{b}_{n/2}^2$ and $\xi_i^2\coloneqq S_i=\hat{a}_1^2+\hat{b}_1^2$ then it follows from (\ref{eqn:deriv_scr_a}) that
\begin{eqnarray}
(X_{2i-1},X_{2i}) & = & (\hat{x}_i,\hat{y}_i) \sqrt{S_i} \sqrt{\frac{1-S_{i-1}/S_i}{S_{n/2}}} \nonumber\\
   & = & (\hat{a}_i,\hat{b}_i) \sqrt{\frac{1-S_{i-1}/S_i}{S_{n/2}}},
\end{eqnarray}
for $i=1,\dots,n/2$ is also a valid random vector in $\mathcal{S}_n$.

Furthermore, since $\sqrt{S_{n/2}}$ constitutes the maximum of $n/2$ variables with linearly increasing distribution, it is $P_{\sqrt{S_{n/2}}}(x)\propto  x^{n-1}$ and it has the support $x\in[0,1)$. This is the distribution of the length of a random vector in the unit ball $\mathcal{B}_{n}$. Since $\sqrt{S_{n/2}}$ is also independent of $\mathbf{X}$ it follows that $\sqrt{S_{n/2}}\,\mathbf{X}$ is a random point in $\mathcal{B}_{n}$.

\section{Runtimes}

It is obvious that a straightforward sorting of the pairs of coordinates has a crucial impact on the performance of our algorithm, since it is of complexity $O(n\log n)$ and dominates the runtime for large $n$. On the other hand, the distribution of the numbers that need to be sorted is known, which can be exploited for optimization. Besides a basic routine that uses the standard \emph{sort()}-function provided for the \CC container \emph{vector} we, therefore, implemented another version using \emph{BucketSort} \cite{bucketsort} leading to an only slightly more complicated code\footnote{Whereas Ref. \cite{bucketsort} prescribes the use of linked lists for the ``buckets'' we here employ arrays (\CC container \emph{vector}) which we find to produce a faster running code.}.
Both algorithms use additional memory, which is why we lastly created an \emph{in situ} version for even $n$. For all methods the required memory was allocated before the measurement. We compare to the method predominantly in use, i.e., Muller's method of rescaling points of a multivariate normal distribution. Since the amount of random numbers required by our method is on average by a factor of $4/\pi$ larger than for Mueller's method, the choice of the random number generator (RNG) has an impact. For this reason we employ several RNGs, namely the \emph{mt19937} Mersenne Twister from the \emph{random}-library in its 32-bit and its 64-bit versions, Marsaglia and Tsang's 64-bit RNG \cite{MARS_RNG}, and the very fast linear congruential RNG \emph{drand48}, whose pseudo random numbers are of lower quality. The values of $n$ at which the measurements were taken are obtained by starting with $n=2$ and either increasing by one if $n$ is even or multiplying with the golden ratio $\phi\approx 1.618034$ and rounding down to the next even value if $n$ is odd.

\begin{figure}
\begin{center}
\includegraphics[width=0.8\columnwidth]{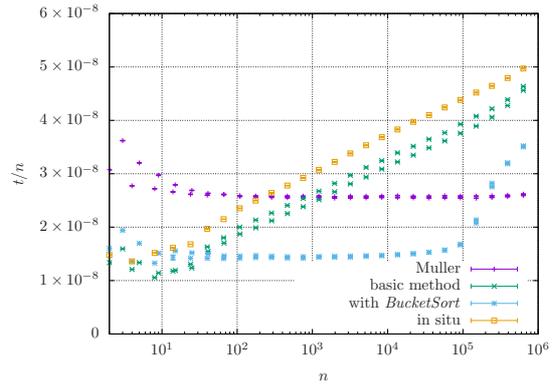}
\end{center}
\caption{\small{\label{fig:times_dr48} Time (in seconds) per component for several algorithms (see text) as function of the number of dimensions $n$ using the \emph{drand48} random number generator. Note the pronounced oscillations between even and odd values of $n$ for the basic version of our method; for small $n$ ($<10$) the ratio $t/n$ is smaller for even $n$ while throughout the interval where the $n\log n$ growth is dominant $t/n$ is smaller if $n$ is odd.}}
\end{figure}

\begin{figure}
\begin{center}
\includegraphics[width=0.8\columnwidth]{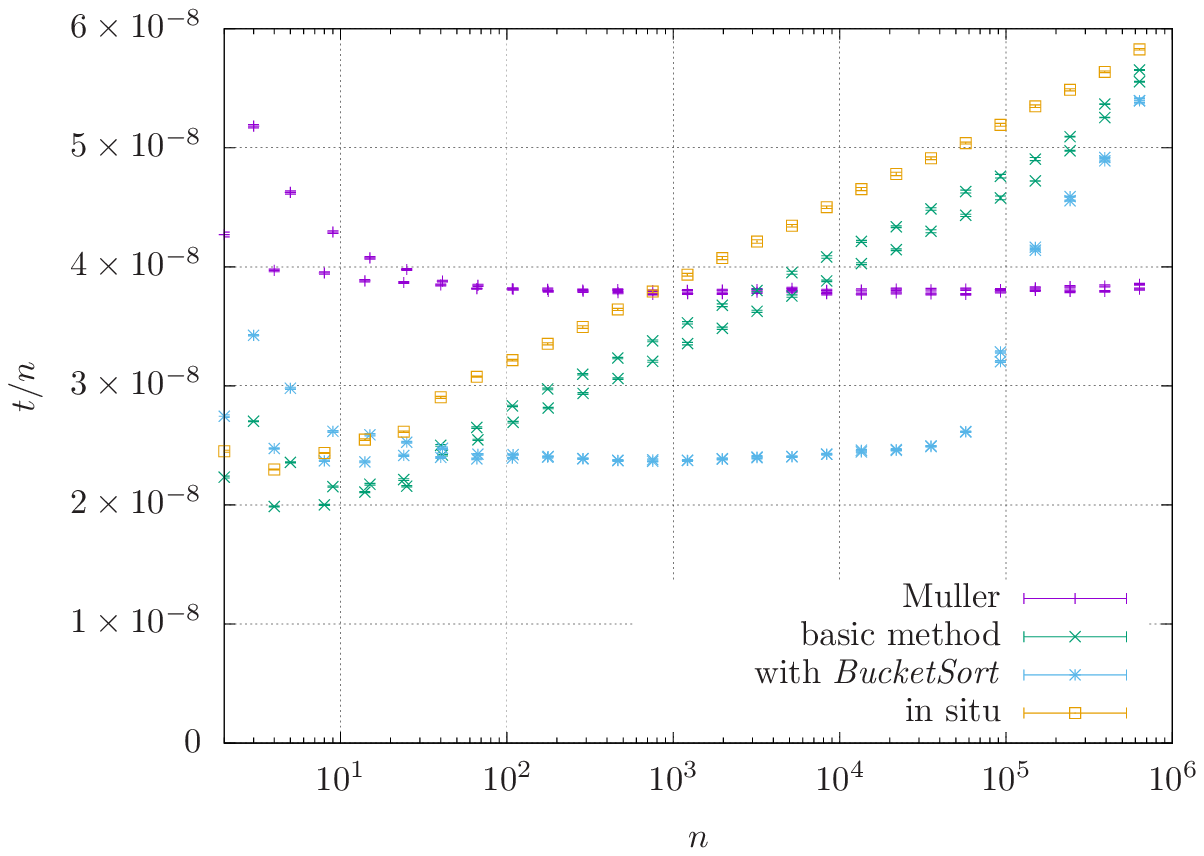}
\end{center}
\caption{\small{\label{fig:times_mt32} Same as Fig.~\ref{fig:times_dr48} but using the 32-bit \emph{Mersenne  Twister} random number generator.}}
\end{figure}

\begin{figure}
\begin{center}
\includegraphics[width=0.8\columnwidth]{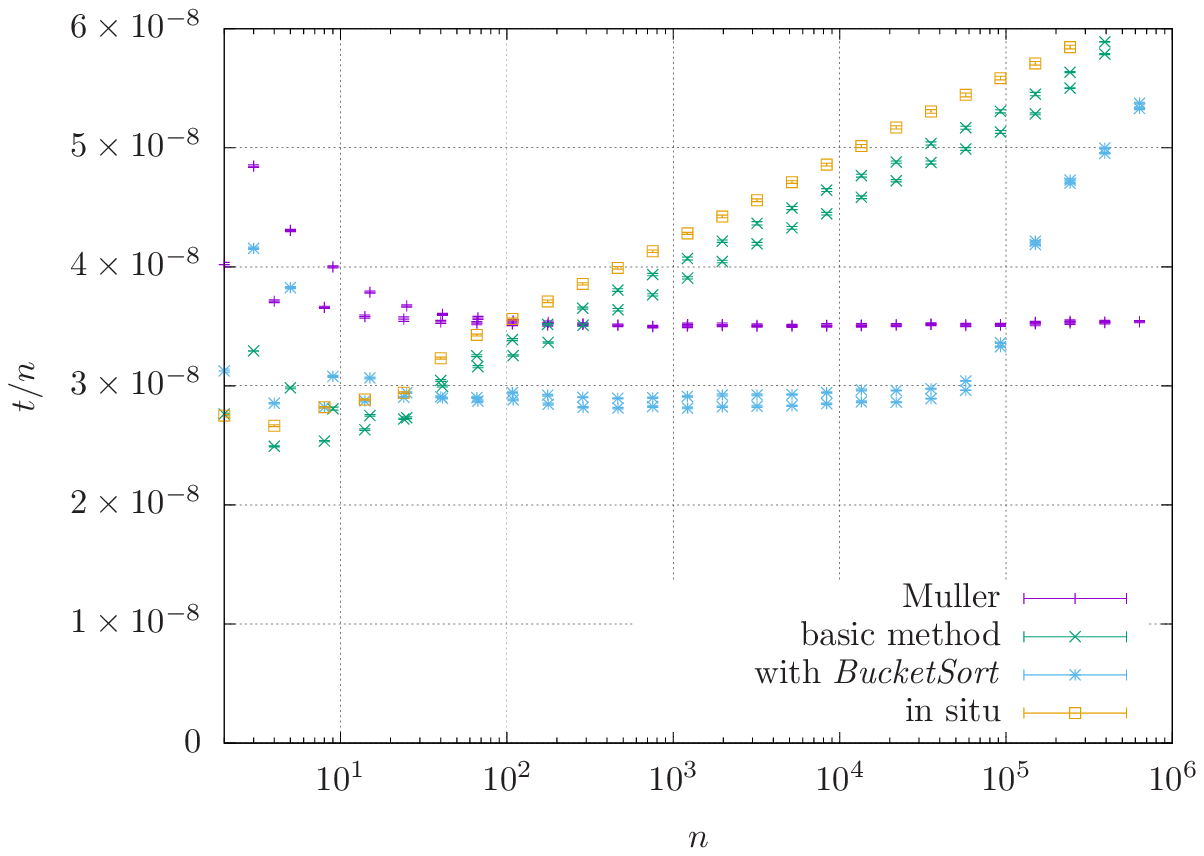}
\end{center}
\caption{\small{\label{fig:times_mt64} Same as Fig.~\ref{fig:times_dr48} but using the 64-bit \emph{Mersenne  Twister} random number generator.}}
\end{figure}

\begin{figure}
\begin{center}
\includegraphics[width=0.8\columnwidth]{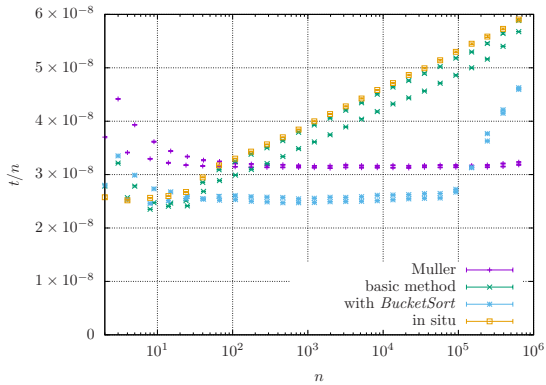}
\end{center}
\caption{\small{\label{fig:times_mars} Same as Fig.~\ref{fig:times_dr48} but using Marsaglia and Tsang's 64-bit random number generator.}}
\end{figure}

The \CC code is available \cite{code_on_gh} on GitHub. It was compiled using the gcc compiler (version 8.3) with O3-optimization and executed on a Intel Core i5-6500 CPU at 3.2GHz. The average times required for obtaining one component of the random vector are displayed in Figs.~\ref{fig:times_dr48}-\ref{fig:times_mars}.

As can be seen, the basic version of our method is faster than Muller's method for small $n$ regardless of the RNG. However, as expected the demands of sorting soon slows down this version and for $n>40$ the usage of \emph{BucketSort} becomes profitable. For all values $n<10^5$ the latter remains faster than Muller's method. We suspect that the slowing down for even larger $n$ is caused by the growing memory demand and an increase in memory latency it entails. Further optimization is most likely possible in this regard. The \emph{in situ} routine is functional very similar to the basic version and it also performs similarly. It is somewhat slower since it requires the calculation of $n$ square roots instead of $n/2$. One odd observation is that Mueller's method runs faster\footnote{When the  code is compiled on our system using optimization, it appears that $\ln(1-{\rm random\_number}())$ and to a lesser extend $\cos\left(2\pi {\rm random\_number}()\right)$ are evaluated faster when the 64-Bit Mersenne Twister is used, even though calling ${\rm random\_number}()$ takes longer in this case.} when combined with the 64-bit Mersenne Twister compared to the 32-bit version. This shows that there are subtle issues related to code design, optimization, choice of compiler, etc. whose discussion would exceed the scope of the paper. We deem it likely that for the implementation of all methods presented here, further optimization is possible.

\section{Conclusion}

It this paper we have suggested a novel method to generate random points on an $n$-sphere that is closely related to but more efficient than Sibuya's method \cite{Sibuya}. It is conceptually very simple and avoids computationally expensive functions. Although in its basic form it has a less favorable algorithmic complexity than the widely used method proposed by Mueller \cite{Muller} this can be mended by the selection of a suitable sorting algorithm and we found it to be faster at least for a number of dimensions smaller than $10^5$, i.e., for most realistic scenarios. 

\section*{Acknowledgement}
This  project  was  funded  by  the  Deutsche  Forschungsgemeinschaft  (DFG,  German  Research  Foundation)  under Project No.\ 189 853 844–SFB/TRR 102 (project B04).

\appendix 
\section{CODE}
\begin{lstlisting}
#include <cmath>
#include <vector>
#include <algorithm>
#include <random>

#ifndef M_PI
#define M_PI 3.14159265358979324
#endif

//std::mt19937 gen;
std::mt19937_64 gen;

inline double random_num(void){
    return static_cast<double>(gen())/gen.max();
}

/*******************Mueller**********************/

void rand_vect_mueller_odd(int n,double* v){
    double r,p,S(0.);
    for(int k=0;k<n-1;k+=2){
        r = sqrt( -2*log(1.0-random_num()) );
        p = random_num();
        v[k] = cos(2*M_PI*p) ;
        v[k+1] = p<0.5 ?  r*sqrt(1.0 -v[k]*v[k]) 
                       : -r*sqrt(1.0 -v[k]*v[k]);
        v[k] *= r;
        S += v[k]  *v[k];
        S += v[k+1]*v[k+1];
    }
    r = sqrt( -2*log(1.0-random_num()) );
    p = 2*M_PI*random_num();
    v[n-1] = r * cos(p) ;
    S = sqrt( S + v[n-1]*v[n-1] );
    for(int k=0;k<n;++k)
        v[k] /= S;
}

void rand_vect_mueller_even(int n,double* v){
    double r,p,S(0.);
    for(int k=0;k<n;k+=2){
        r = sqrt( -2*log(1.0-random_num()) );
        p = random_num();
        v[k] = cos(2*M_PI*p) ;
        v[k+1] = p<0.5 ?  r*sqrt(1.0 -v[k]*v[k]) 
                       : -r*sqrt(1.0 -v[k]*v[k]);
        v[k] *= r;
        S += v[k]  *v[k];
        S += v[k+1]*v[k+1];
    }
    S = sqrt(S);
    for(int k=0;k<n;++k)
        v[k] /= S;
}

void rand_vect_mueller(int n,double* v){
    return n&1 ? rand_vect_mueller_odd(n,v)
               : rand_vect_mueller_even(n,v);
}

/*******************Basic************************/

struct triple{
    double S,a,b;
    bool operator<(const triple& other) const{
        return (S < other.S);   }
};

void rand_vect_basic_even(int n,double* v){
    static std::vector< triple > twoDpts;
    if( twoDpts.size() != n/2 )
        twoDpts.resize(n/2 , triple() );
    for(int i=0;i<n/2;++i){
        do{
            twoDpts[i].a = 2*random_num()-1.0;
            twoDpts[i].b = 2*random_num()-1.0;
            twoDpts[i].S = twoDpts[i].a*twoDpts[i].a
                         + twoDpts[i].b*twoDpts[i].b;
        }
        while( twoDpts[i].S > 1.0 );
    }
    std::sort( twoDpts.begin() , twoDpts.end() );
    double t;
    double s = 1.0 / twoDpts[n/2-1].S ;
    double S_I,S_II(0.0);
    for(int i=0;i<n/2;++i){
        S_I  = S_II ;
        S_II = twoDpts[i].S ;
        t = sqrt( (1.0 - S_I/S_II )*s ) ;
        v[2*i]   = twoDpts[i].a * t ;
        v[2*i+1] = twoDpts[i].b * t ;
    }
    return;
}

void rand_vect_basic_odd(int n,double* v){
    ++n;
    static std::vector< triple > twoDpts;
    if( twoDpts.size() != n/2 )
        twoDpts.resize( n/2 , triple() );
    for(int i=0;i<n/2;++i){
        do{
            twoDpts[i].a = 2*random_num()-1.0;
            twoDpts[i].b = 2*random_num()-1.0;
            twoDpts[i].S = twoDpts[i].a*twoDpts[i].a
                         + twoDpts[i].b*twoDpts[i].b;
        }
        while( twoDpts[i].S > 1.0 );
    }
    std::sort( twoDpts.begin() , twoDpts.end() );
    double t, s( 1.0 / twoDpts[n/2-1].S );
    s = s / ( 1.0 - s * twoDpts[0].a * twoDpts[0].a );
    double S_I(0.0), S_II(twoDpts[0].S);
    v[0] = twoDpts[0].b * sqrt( s ) ;
    for(int i=1;i<n/2;++i){
        S_I  = S_II ;
        S_II = twoDpts[i].S ;
        t = sqrt( (1.0 - S_I/S_II )*s ) ;
        v[2*i-1] = twoDpts[i].a * t ;
        v[2*i]   = twoDpts[i].b * t ;
    }
    return;
}

void rand_vect_basic(int n,double* v){
    return n&1 ? rand_vect_basic_odd(n,v)
               : rand_vect_basic_even(n,v);
}

/*******************with BucketSort**************/

const int BUCKET_SIZE=16;

void rand_vect_bs_even(int n,double* v){
    static std::vector<std::vector<triple> > twoDpts;
    int nb = 1+n/BUCKET_SIZE; //
    twoDpts.resize( nb , std::vector<triple>(0) );
    triple tr;
    for(int i=0;i<n/2;++i){
        do{
            tr.a = 2*random_num()-1.0;
            tr.b = 2*random_num()-1.0;
            tr.S = tr.a*tr.a + tr.b*tr.b ;
        }
        while( tr.S > 1.0 );
        twoDpts[tr.S*nb].push_back( tr );
    }
    for (int bi = 0; bi < nb; ++bi)
        sort(twoDpts[bi].begin(),twoDpts[bi].end());
    
    int last_bucket = nb-1;
    while(!twoDpts[last_bucket ].size())
        --last_bucket;
    double t;
    double S_I,S_II(0.0);
    double s = 1.0 / twoDpts[last_bucket].rbegin()->S;
    int i=0;
    for (int bi = 0; bi <=last_bucket; ++bi){
        for (int j = 0; j<twoDpts[bi].size(); ++j){
            S_I  = S_II ;
            S_II = twoDpts[bi][j].S ;
            t = sqrt( (1.0 - S_I/S_II )*s ) ;
            v[i++] = twoDpts[bi][j].a * t ;
            v[i++] = twoDpts[bi][j].b * t ;
        }
        twoDpts[bi].resize(0);
    }
    return;
}

void rand_vect_bs_odd(int n,double* v){
    ++n;
    static std::vector< std::vector<triple>  > twoDpts;
    int nb = 1+n/BUCKET_SIZE; //
    twoDpts.resize( nb , std::vector<triple>(0) );
    triple tr;
    for(int i=0;i<n/2;++i){
        do{
            tr.a = 2*random_num()-1.0;
            tr.b = 2*random_num()-1.0;
            tr.S = tr.a*tr.a + tr.b*tr.b ;
        }
        while( tr.S > 1.0 );
        twoDpts[tr.S*nb].push_back( tr );
    }
    for (int bi = 0; bi < nb; ++bi)
        sort(twoDpts[bi].begin(),twoDpts[bi].end());
    
    int first_bucket = 0;
    while(!twoDpts[first_bucket ].size()) 
        ++first_bucket;
    int last_bucket = nb-1;
    while( !twoDpts[last_bucket].size() )
        --last_bucket;
    double t;
    double s = 1.0 / twoDpts[last_bucket].rbegin()->S ;
    s = s/( 1.0 - s*twoDpts[first_bucket].begin()->a
                   *twoDpts[first_bucket].begin()->a);
    double S_I,S_II;
    v[0] = twoDpts[first_bucket].begin()->b * sqrt( s ) ;
    S_II = twoDpts[first_bucket].begin()->S ;
    int i=1;
    int j=1;
    for (int bi = first_bucket; bi <=last_bucket; ++bi){
        for ( ;j<twoDpts[bi].size(); ++j){
            S_I  = S_II ;
            S_II = twoDpts[bi][j].S ;
            t = sqrt( (1.0 - S_I/S_II )*s ) ;
            v[i++] = twoDpts[bi][j].a * t ;
            v[i++] = twoDpts[bi][j].b * t ;
        }
        j=0;
        twoDpts[bi].resize(0);
    }
    return;
}

void rand_vect_bs(int n,double* a){
    return n&1 ? rand_vect_bs_odd(n,a)
               : rand_vect_bs_even(n,a);
}

/*******************inSitu***********************/

struct my_pair{
    double S,b;
    bool operator<(const my_pair& other) const{
        return ( fabs(S) < fabs(other.S) ); }
};

void rand_vect_in_situ_even(int n,double* v){
    my_pair* twoDpts( reinterpret_cast<my_pair*>( v ) );
    for(int i=0;i<n/2;++i){
        do{
            twoDpts[i].S = 2*random_num()-1.0;
            twoDpts[i].b = 2*random_num()-1.0;
            // store squared radius and sign of a in S:
            twoDpts[i].S = (twoDpts[i].S>0)
                ?   twoDpts[i].S*twoDpts[i].S
                    + twoDpts[i].b*twoDpts[i].b
                : -(twoDpts[i].S*twoDpts[i].S
                    + twoDpts[i].b*twoDpts[i].b) ;
        }
        while( fabs(twoDpts[i].S) > 1.0 );
    }
    std::sort( twoDpts , twoDpts + n/2 );
    double t;
    double s = 1.0 / fabs(twoDpts[n/2-1].S) ;
    double S_II(0.0);
    for(int i=0;i<n/2;++i){
        t  = S_II ;
        S_II =fabs( twoDpts[i].S) ;
        t = sqrt( (1.0 - t/S_II )*s ) ;
        twoDpts[i].S = twoDpts[i].S>0
            ?  sqrt(  twoDpts[i].S
                - twoDpts[i].b*twoDpts[i].b ) * t
            : -sqrt( -twoDpts[i].S
                - twoDpts[i].b*twoDpts[i].b ) * t;
        twoDpts[i].b = twoDpts[i].b * t ;
    }
    return;
}

/************************************************/

\end{lstlisting}

\end{document}